\documentclass[10pt]{article}
\usepackage{amsthm,amsmath,amsfonts}
\usepackage{fullpage,enumerate}
\usepackage{tikz}
\usepackage[colorlinks,citecolor=blue,urlcolor=blue]{hyperref}
\usepackage{hypernat}
\usepackage[notref,notcite]{showkeys}
\usepackage{graphicx}

\begin{document}

\numberwithin{equation}{section}

\title{The waiting time for a second mutation: an alternative to the Moran model}

\bigskip
\author{  Rinaldo B. Schinazi\\
{\it University of Colorado, Colorado Springs}}
\maketitle

email: rschinaz@uccs.edu

fax: 1 719 2553605

\medskip

{\bf Abstract} The appearance of cancer in a tissue is thought to be the result of two or more successive mutations. We propose a stochastic model that allows for an exact computation of the distribution of the waiting time for a second mutation. This models the time of appearance of the first cancerous cell in a tissue. Our model is an alternative to the Moran model with mutations.

\section {The model}

The idea of successive mutations to trigger the appearance of a cancerous cell goes back to at least Muller (1951). The first mathematical model proposed for this phenomenon goes back to Armitage and Doll (1954). There has been a great deal of discussion on the number of successive mutations necessary to get a cancerous cell. It seems that this number depends on the organ, see Knudson (1971) and Moolgavkar and Luebeck (1992) . However, some authors have argued that two mutations models are flexible enough to model most cancers, see Armitage and Doll (1957) and Moolgavkar and Knudson (1981). This is the point of view we adopt.

We now describe our model.
We are interested in the time it takes for a given organ to have a first cancerous cell. We assume that all cells are in one of three stages: healthy, pre-cancerous (i.e. type 1) and cancerous (i.e. type 2). We start the process with all cells healthy. As the cells divide pre-cancerous cells may appear due to a type 1 mutation on a healthy cell. A type 2 mutation on a pre-cancerous cell makes the cell cancerous. 

The number of type 1 mutations is modeled by a Poisson process with rate $\mu_1 N$. We think of $\mu_1$ as a mutation rate and $N$ as a division rate. Every time a type 1 mutation appears there is a chance that a type 2 mutation appears.  We model the appearance of a type 2 mutation by using exponential random variables with rate $\mu_2$. 
More precisely, let $N_1(t)$ be the number of type 1 mutations that have occurred up to time $t$. Let $T_1<T_2<\dots$ be the arrival times of this Poisson process. Given that $N_1(t)=k$ and that $T_1=t_1$, $T_2=t_2, \dots ,T_k=t_k$ let  $S_1,S_2,\dots, S_k$ be random variables with density 
$$P(S_i>t|T_i=t_i)=\exp(-\mu_2(t-t_i))\mbox{ for }t>t_i,$$
and $i=1,\dots,k$. We also assume that given $T_1=t_1$, $T_2=t_2\dots T_k=t_k$, the random variables 
$S_1,S_2,\dots, S_k$ are independent. The random variables $S_i$ are the times when a type 2 mutation appears. The minimum of these times (i.e. the first time a type 2 mutation appears in the tissue) is denoted by $\tau_2$.

Here is our main result.

\medskip

{\bf Theorem 1. }{\sl Let $\tau_2$ be the time for the first type 2 mutation to appear. Then,
$$P(\tau_2>t)=\exp[\mu_1 Nt(-1+\frac{1-\exp(-t\mu_2)}{t\mu_2})].$$}

\medskip

Note that 
$$P(\tau_2>t)=\exp[\mu_1 tf(t\mu_2)]^N$$
where
$f(x)=-1+\frac{1-\exp(-x)}{x}$.

This shows that our model is equivalent to a model with $N$ independent cells. Observe also that the function $f$ decreases from 0 to $-1$ for $x$ in $[0,+\infty)$. In particular, $\mu_1$ and $\mu_2$ do not hold symmetric roles in the formula. For instance, a small $\mu_1$   cannot be compensated by a large $\mu_2$.

\section{Discussion}

Waiting times for successive mutations have been recently studied by several authors using the Moran model with mutations, see Iwasa et al. (2004) and (2005) and Wodarz and Komarova (2005). A more mathematical approach is taken by Durrett et al. (2009) and Schweinsberg (2008). The Moran model assumes a fixed number $N$ of cells. Each cell lives for a mean 1 exponential time and then is replaced by a new cell chosen at random from one of the $N$ cells. Moreover, a healthy cell mutates into a precancerous cell at rate $\mu_1$ and a precancerous mutates into a cancerous cell at rate $\mu_2$. Each new cell has the same number of mutations as its parent. An attractive feature of the Moran model is that it is defined as a cell based model. A drawback is that the analysis of the model is quite involved and the only results that can be hoped for are non rigorous approximations or limits for different configurations of $N$, $\mu_1$ and $\mu_2$. In contrast, our model gives an exact formula for the waiting time distribution.

Our model is not fundamentally different from the Moran model. In particular, in the Moran model a precancerous cell has equal death and birth rates. Hence, a precancerous cell generates a process which in average does not gain or lose cells. It turns out that estimating the stochastic fluctuations around this expected value is quite difficult.  In our model we assume instead that a precancerous cell does not die or give birth, it stays put waiting for a type 2 mutation. Therefore, the expected offspring of a precancerous cell is also 1 in our model but there are no stochastic fluctuations to estimate.

As a consequence of Theorem 1 we have the following limits. To compute these limits assume that $\mu_1$ and $\mu_2$ are functions of $N$.

$\bullet$ Assume that 
$$\lim_{N\to\infty}\frac{\mu_2} {\mu_1 N}=\alpha\in (0,+\infty).$$
Then,
$$\lim_{N\to\infty}P(\mu_1 N\tau_2>t)=\exp[t(-1+\frac{1-\exp(-t\alpha)}{t\alpha})].$$

\medskip

$\bullet$ Assume that 
$$\lim_{N\to\infty}\frac{\mu_2} {\mu_1 N}=+\infty.$$
Then,
$$\lim_{N\to\infty}P(\mu_1 N\tau_2>t)=\exp(-t).$$

\medskip

$\bullet$ Assume that 
$$\lim_{N\to\infty}\mu_1 N\mu_2=\alpha\in (0,+\infty)\mbox{ and }\lim_{N\to\infty}\mu_2=0.$$
Then,
$$\lim_{N\to\infty}P(\tau_2>t)=\exp(-\frac{1}{2}\alpha t^2).$$

\medskip

In particular, the distribution of $\tau_2$ exhibits at least three different behaviors depending on the relative magnitude of $\mu_1$, $\mu_2$ and $N$. These limits show a number of similarities with limits found for the Moran model by Durrett, Schmidt and Schweinsberg (2009) and Schweinsberg (2008).

Note also that as $t$ approaches 0 (for fixed $\mu_1$, $\mu_2$ and $N$) 
$$P(\tau_2\leq t)\sim \frac{1}{2}\mu_1\mu_2 N t^2.$$
This is consistent with the model of Armitage and Doll (1954).

\section{ The proof of Theorem 1}

Let $N_1(t)$ be the number of type 1 mutations that occurred in the tissue up to time $t$. Given $\{ N_1(t)=k\}$ no type 2 mutation has occurred (i.e. $\{\tau_2>t\}$) if and only if none of the $k$ type 2 mutation exponential random variables have occurred. Denote these $k$ random variables by $S_1,S_2,\dots, S_k$. We have
$$P(\tau_2>t|N_1(t)=k)=P(S_1>t,S_2>t,\dots, S_k>t|N_1(t)=k).$$
 Let $T_1,T_2,\dots,T_k$ be the arrival times of the Poisson process $N_1$. By definition of the random variables $S_1,S_2,\dots, S_k$ we have for $i=1,\dots,k$
$$P(S_i>t|T_i=t_i)=\exp(-\mu_2(t-t_i))\mbox{ for }t>t_i.$$
Given $(T_1,T_2,\dots,T_k)$ the random variables $S_1,S_2,\dots, S_k$ are  conditionally independent. Hence, 
$$P(\tau_2>t|N_1(t)=k)=\int_{0<t_1<t_2<\dots<t_k<t}\exp(-\mu_2(t-t_1))\dots\exp(-\mu_2(t-t_k))f(t_1,t_2,\dots,t_k)dt_1dt_2\dots dt_k,$$
where $f$ is the density of the random vector $(T_1,T_2,\dots,T_k)$ conditioned on $\{ N_1(t)=k\}$. A classical Poisson process result is that this conditional distribution is the order statistics distribution corresponding to $k$ independent random variables uniformly distributed on $(0,t)$, see for instance Proposition 5.6 in Bhattacharya and Waymire (1990). Therefore,

$$P(\tau_2>t|N_1(t)=k)=\exp(-\mu_2 k t)\int_{0<t_1<t_2<\dots<t_k<t}\frac{k!}{t^k}\exp(\mu_2(t_1+t_2+\dots+t_k)dt_1\dots dt_k.$$
In order to compute this integral we make the following remark.
Let $U_1, U_2,\dots,U_k$ be independent and uniformly distributed on $(0,t)$. Let $U_{(1)}<U_{(2)}<\dots<U_{(k)}$ be the corresponding order statistics. We have that
$$E[\exp(\mu_2(U_{(1)}+U_{(2)}+\dots+U_{(k)}))]=\int_{0<u_1<u_2<\dots<u_k<t}\frac{k!}{t^k}\exp(\mu_2(u_1+u_2+\dots+u_k)du_1\dots du_k.$$
Observe that
$$U_{(1)}+U_{(2)}+\dots+U_{(k)}=U_1+ U_2+\dots+U_k.$$
Hence,
$$E[\exp(\mu_2(U_{(1)}+U_{(2)}+\dots+U_{(k)}))]=E[\exp(\mu_2(U_1+ U_2+\dots+U_k))]=E[\exp(\mu_2 U_1)]^k.$$
It is easy to compute
$$E[\exp(\mu_2 U_1)]=\frac{1}{t\mu_2}(\exp(t\mu_2)-1).$$
Therefore,
$$P(\tau_2>t|N_1(t)=k)=\exp(-\mu_2 k t)[\frac{1}{t\mu_2}(\exp(t\mu_2)-1)]^k=[\frac{1}{t\mu_2}(1-\exp(-t\mu_2))]^k.$$

Now,
$$P(\tau_2>t)=\sum_{k=0}^\infty P(N_1(t)=k)P(\tau_2>t|N_1(t)=k)=\sum_{k=0}^\infty \exp(-\mu_1 Nt)\frac{(\mu_1 Nt)^k}{k!}[\frac{1}{t\mu_2}(1-\exp(-t\mu_2))]^k.$$
Summing the series yields
$$P(\tau_2>t)=\exp(-\mu_1 Nt)\exp[\frac{\mu_1 N}{\mu_2}(1-\exp(-t\mu_2))].$$
This formula can be rewritten as

$$P(\tau_2>t)=\exp[\mu_1 Nt(-1+\frac{1-\exp(-t\mu_2)}{t\mu_2})].$$

The proof of Theorem 1 is complete.

\bigbreak

{\bf References}

\medskip

P. Armitage P. and R. Doll (1954) The age distribution of cancer and a multistage theory of carcinogenesis. British Journal of cancer {\bf 8:} 1-12.

P. Armitage P. and R. Doll (1957) A two-stage theory of carcinogenesis in relation to the age distribution of human cancer. British Journal of cancer {\bf 11:} 161-169.

R.N.Bhattacharya and E.C.Waymire (1990){\it Stochastic processes with applications}, Wiley.

R. Durrett, D. Schmidt, and J. Schweinsberg (2009). A waiting time problem arising from 
the study of multi-stage carcinogenesis. Annals of Applied probability {\bf 19}, 676-718.

H.W. Hethcote  and A.G. Knudson Jr (1978). Model for the incidence of embryonal cancers: Application to retinoblastoma. Proceedings of the National Academy of Sciences {\bf 75}, 2453-2457.

Y. Iwasa, F. Michor, N.L. Komarova, and M.A. Nowak (2005). Population 
genetics of tumor suppressor genes. J. Theoret. Biol. {\bf 233}, 15-23

Y. Iwasa, F. Michor, and M.A. Nowak (2004). Stochastic tunnels in evolutionary 
dynamics. Genetics {\bf 166}, 1571-1579. 

A.G.Knudson (2001) Two genetic hits (more or less) to cancer. Nature reviews Cancer {\bf 1: }157-162.

A. G. Knudson (1971). Mutation and cancer: statistical study of retinoblastoma. Proc. Natl. 
Acad. Sci. USA {\bf 68}, 820-823. 

N.L. Komarova, A. Sengupta,  and M.A. Nowak (2003). Mutation-selection 
networks of cancer initiation: Tumor suppressor genes and chromosomal insta- 
bility. J. Theoret. Biol. {\bf 223}, 433-450.

S.H. Moolgavkar and A.G. Knudson (1981). Mutation and cancer:a model for human carcinogenesis. J.Natl. Cancer. Inst {\bf 66}, 1037-1052.

S. H. Moolgavkar and E. G. Luebeck (1992). Multistage carcinogenesis: population-based 
model for colon cancer. J. Natl. Cancer Inst. {\bf 18}, 610-618. 

J. Schweinsberg (2008). The waiting time for m mutations. Electron. J. Probab. 
{\bf 13}, 1442-1478.

D. Wodarz and N.L. Komarova (2005). {\it Computational Biology Of Cancer.} Lecture Notes And Mathematical Modeling. World ScientiÞc, Singapore.

\end{document}